\documentclass[a4paper,12pt]{article}
\usepackage[french, ngerman, italian, english]{babel}
\usepackage{amsmath,amsthm, amssymb}
\usepackage{setspace}
\usepackage{anysize}
\usepackage{url}
\usepackage{graphicx} 
\usepackage[colorlinks=true]{hyperref}
\usepackage{hyperref}
\usepackage{mathptmx} 
\usepackage{paralist}
\usepackage{verbatim}
\usepackage[utf8]{inputenc}

\theoremstyle{plain}
\newtheorem{thm}{\bf Theorem}[section] 

\newtheorem{proposition}[thm]{\bf Proposition}
\newtheorem{thmnonumber}{\bf Main Theorem}

\newtheorem{cor}[thm]{\bf Corollary}
\newtheorem*{cor*}{\bf Corollary}
\theoremstyle{definition}

\theoremstyle{remark}
\newtheorem{remark}[thm]{\bf Remark}
\newtheorem{example}[thm]{\bf Example}

\begin{document}
\title{Collapses, products and LC manifolds}
\author{Bruno Benedetti\thanks{%
Supported by DFG via the Berlin Mathematical School} \\
\small Inst.\ Mathematics, MA 6-2, TU Berlin\\
\small 10623 Berlin, Germany\\
\small \url{benedetti@math.tu-berlin.de} 
} 
\date{May 16, 2010}
\maketitle

\begin{abstract}\noindent 
The quantum physicists Durhuus and Jonsson (1995) introduced the class of ``locally constructible'' (LC) triangulated manifolds and showed that all the LC $2$- and $3$-manifolds are spheres. 
We show here that for each $d > 3$ some LC $d$-manifolds are not spheres. We prove this result by studying how to collapse products of manifolds with one facet removed.
\end{abstract}

\section{Introduction}
\emph{Collapses} are a classical notion in Combinatorial Topology, 
originally introduced in the Thirties by Whitehead \cite{Whitehead}, extensively studied in the Sixties by Bing, Cohen, Lickorish and Zeeman among others,
yet also at the center of recent works such as \cite{BarmakMinian} and \cite{Crowley}. 

Given a polytopal (or a regular CW) complex, a collapse is a move that cancels two faces and yields a smaller complex which is topologically a strong deformation retract of the starting one. Any complex that is collapsible (i.e. transformable into a point via a sequence of collapses) is thus also contractible. Conversely, every shellable contractible complex is collapsible.

However, not all contractible complexes are collapsible: A famous two-dimensional counterexample is given by Zeeman's dunce hat \cite{ZeemanDunceHat}.
According to the work of Whitehead \cite{Whitehead} and Cohen \cite{Cohen1}, a complex $C$ is contractible if and only if some collapsible complex $D$ collapses also onto $C$. In fact, one can construct a collapsible triangulated $3$-ball with only $8$ vertices that collapses onto a copy of the dunce hat \cite{BenedettiLutz}. 
Cohen's result is obtained by taking \emph{products}: Zeeman \cite{ZeemanDunceHat} first noticed that the product of the dunce hat with a segment $\mathbb{I}$ is polyhedrally collapsible and asked whether the same holds for any contractible $2$-complex. 
(The question, known as \emph{Zeeman's conjecture}, is still open \cite{HAMS}. For polyhedral collapsibility, see e.g.~\cite[pp.~42--48]{Hudson}.) 

Cohen \cite[Corollaries~3~\&~4]{Cohen1} showed that the product of any contractible $d$-complex $C$ with the $q$-di\-mensional cube $\mathbb{I}^q$ polyhedrally collapses onto a point, provided $q \ge  \max(2d, 5)$. At the same time, $C \times \mathbb{I}^q$ collapses onto $C$ (cf.~Corollary \ref{cor:cohen}).

\newpage
It was first discovered by Bing \cite{BING} that some triangulations of $3$-balls are not collapsible. For each $d \ge 3$, Lickorish \cite{LICK} proved that some triangulated $d$-balls of the form $S - \Delta$ (with $S$ a $d$-sphere and  $\Delta$ a facet of $S$) are not collapsible. Bing's and Lickorish's claim were recently strengthened by the author and Ziegler \cite[Thm.~2.19]{BZ}, who showed that for each $d \ge 3$ certain triangulated $d$-balls of the form $S - \Delta$ do not even collapse onto any $(d-2)$-dimensional subcomplex of $S$. 
These three results were all obtained via knot theory. In fact,  a $3$-ball may contain arbitrarily complicated three-edge-knots in its $1$-skeleton. Depending on how complicated the knot is, one can draw sharp conclusions on the collapsibility of the $3$-ball and of its successive suspensions.

In the Nineties, two quantum physicists, Durhuus and Jonsson \cite{DJ}, introduced  the term ``\emph{LC $d$-manifold}'' to describe a manifold that can be obtained from a tree of $d$-polytopes by repeatedly identifying two combinatorially equivalent adjacent $(d-1)$-faces in the boundary ($d \ge 2$).  Plenty of spheres satisfy this bizarre requirement: In fact, all shellable and all constructible $d$-spheres are LC (cf.~\cite{BZ}). At the same time, simplicial LC $d$-manifolds are only exponentially many when counted with respect to the number of facets, while arbitrary (simplicial) $d$-manifolds are much more numerous \cite[Chapter~2]{Benedetti-diss}.

Durhuus and Jonsson noticed that the class of LC $d$-manifolds coincides with the class of all $d$-spheres for $d=2$. But what about higher dimensions?

For $d=3$, they were able to prove one of the two inclusions, namely, that all LC $3$-manifolds are spheres \cite[Theorem~2]{DJ}. The other inclusion does not hold: For each $d \ge 3$, some $d$-spheres are not LC, as established in \cite{BZ}. The examples of non-LC spheres are given by $3$-spheres with a three-edge-knot in their $1$-skeleton (provided the knot is sufficiently complicated!) and by their successive suspensions.

The analogy with the aforementioned obstructions to collapsibility is not a coincidence: In fact, the LC $d$-spheres can be characterized \cite[Theorem~2.1]{BZ} as the $d$-spheres that collapse onto a $(d-2)$-complex after the removal of a facet. (It does not matter which facet you choose.) This characterization can be easily extended to (closed) manifolds:
\bigskip

\emph{
A $d$-manifold $M$ is LC if and only if $M$ minus a facet collapses onto a $(d-2)$-complex.}

\bigskip

\noindent Exploiting this characterization, in the present paper we prove the following statement: 

\begin{thmnonumber}
\label{mainthm:1}
The product of LC manifolds is an LC manifold.
\end{thmnonumber}

\noindent The proof, which is elementary, can be outlined as follows: Suppose
a manifold $M$ (resp.~$M'$) minus a facet collapses onto a $(\dim \, M-2)$-complex $C$ (resp.~a $(\dim \, M'-2)$-complex $C'$). We show that the complex obtained by removing a facet from $M \times M'$ collapses onto the complex $(C \times M') \,\cup \; (M \times C')$, which is $(\dim M + \dim M' -2)$-dimensional. 
 
 As a corollary, we immediately obtain that some LC $4$-manifolds are not spheres, but rather products of two LC $2$-spheres. This enables us to solve Durhuus--Jonsson's problem for all dimensions:

\begin{thmnonumber} \label{mainthm:2}
The class of LC $2$-manifolds coincides with the class of all $2$-spheres. 

\noindent The class of LC $3$-manifolds is strictly contained in the class of all $3$-spheres. 
 
\noindent For each $d \ge 4$, the class of LC $d$-manifolds and the class of all $d$-spheres are overlapping, but none of them is contained in the other.
\end{thmnonumber}

By the work of Zeeman (see e.g. \cite{BJOE}), for every positive integer $d$, every shellable or constructible $d$-manifold is a $d$-sphere. Thus, the properties of  shellability and constructibility are obviously \emph{not} inherited by products. All $2$-spheres are LC, constructible and shellable; however, for each $d \ge 3$, all shellable $d$-spheres are constructible, all constructible $d$-spheres are LC, but some LC $d$-spheres are not constructible \cite{BZ}. It is still unknown whether all constructible spheres are shellable.

\subsection{Definitions}
A \emph{polytopal complex} is a finite, nonempty collection $C$ of
polytopes (called the \emph{faces} of $C$) in some Euclidean space
$\mathbb{R}^k$, such that (1) if $\sigma$ is a polytope in $C$ then all the faces of $\sigma$ are elements of $C$ and (2) the intersection of any two polytopes of $C$ is a face of both.
If $d$ is the largest dimension of a polytope of $C$, the polytopal complex $C$ is called \emph{$d$-complex}. An inclusion-maximal face of $C$ is called \emph{facet}. 
A $d$-complex is \emph{simplicial} (resp. \emph{cubical}) if all of its facets are simplices (resp. cubes). 
Given an $a$-complex $A$ and a $b$-complex $B$, the \emph{product} $C = A \times B$ is an $(a+b)$-complex whose nonempty faces are the products $P_\alpha \times P_\beta$, where $P_\alpha$ (resp. $P_\beta$) ranges over the nonempty polytopes of $A$ (resp. $B$). In general, the product of two simplicial complexes is \emph{not} a simplicial complex, while the product of two cubical complexes is a cubical complex.

Let $C$ be a $d$-complex. An \emph{elementary collapse} is the simultaneous removal from $C$ of a pair of faces $(\sigma, \Sigma)$, such that $\sigma$ is a proper face of $\Sigma$ and of no other face of $C$. (This is usually abbreviated as ``$\sigma$ is a free face of $\Sigma$''; some complexes have no free faces.)
We say the complex $C$ \emph{collapses onto} the complex $D$, and write $C \searrow D$, if $C$ can be deformed onto $D$ by a finite (nonempty) sequence of elementary collapses. Without loss of generality, we may assume that in this sequence the pairs 
$( (d-1)\textrm{-face} , \; d\textrm{-face})$ are removed first; we may also assume that after that, the pairs $( (d-2)\textrm{-face} , \; (d-1)\textrm{-face})$ are removed; and so on. 
A \emph{collapsible} $d$-complex is a $d$-complex that can be collapsed onto a 
single vertex. If $C$ collapses onto $D$, then $D$ is a strong deformation retract of $C$, so $C$ and $D$ have the same homotopy type. In particular, all collapsible
complexes are contractible.

The \emph{underlying space} $|C|$ of a $d$-complex $C$ is the union of all of its faces.
A $d$-\emph{sphere} is a $d$-complex whose underlying space is homeomorphic to 
$
\{ \mathbf{x} \in \mathbb{R}^{d+1} :
   \mathbf{|x|} = 1 \} 
$.
A $d$-\emph{ball} is a $d$-complex with underlying space homeomorphic to 
$
\{ \mathbf{x} \in \mathbb{R}^{d} :
   \mathbf{|x|} \le 1 \} 
$; a \emph{tree of $d$-polytopes} is a $d$-ball whose dual graph is a tree.  
With abuse of language, by \emph{$d$-manifold} we will mean any $d$-complex whose underlying space is homeomorphic to a compact connected topological manifold (without boundary). 

A \emph{locally constructible} (LC) $d$-manifold is a $d$-manifold
obtained from a tree of polytopes by repeatedly identifying a pair of adjacent $(d-1)$-faces of the boundary. (``Adjacent'' means here ``sharing at least a $(d-2)$-face'' and represents a dynamic requirement: after each identification, new pairs of boundary facets might become adjacent and may be glued together.) Equivalently \cite[Theorem~2.1]{BZ} {\cite[Thm.~5.2.6]{Benedetti-diss}}, an LC $d$-manifold is a $d$-manifold that after the removal of a facet collapses onto a $(d-2)$-dimensional subcomplex. For the definition of shellability or constructibility, see e.g.~Bj\"{o}rner \cite[p.~1854]{BJOE}.

\section{Proof of the main results}
In this section, we exploit the characterization of LC manifolds mentioned in the Introduction to prove Main Theorems \ref{mainthm:1} and \ref{mainthm:2}. In fact:
\begin{compactitem}[ -- ]
\item Main Theorem \ref{mainthm:1}
will be a straightforward consequence of Corollary \ref{cor:collapseproduct};
\item  Main Theorem \ref{mainthm:2} follows directly from Remark \ref{rem:perelman}, because we already know that all LC $2$- and $3$-manifolds are spheres \cite[Theorem~2]{DJ}, that all $2$-spheres are LC \cite{DJ} and that some $d$-spheres are not LC for each $d \ge 3$ \cite{BZ}.
\end{compactitem}

\noindent Let us start with a classical result on collapses and products:

\begin{proposition}[Cohen {\cite[p.~254]{Cohen1}}, see also Welker {\cite[Theorem 2.6]{Welker}}] \label{prop:collapseproduct0} Let $A$ and $B$ be two polytopal complexes.
If $A$ collapses onto a complex $C_A$ then $A \times B$ collapses onto $C_A \times B$. 
\end{proposition}

\begin{proof}
Let $B_1, \ldots, B_M$ be an ordered list of all the faces of $B$, ordered by weakly decreasing dimension. Let $(\sigma_1^A, \Sigma_1^A)$ be the first pair of faces appearing in the collapse of $A$ onto $C_A$. We perform the $M$ collapses
$(\sigma_1^A \times B_1, \Sigma_1^A \times B_1)$, 
$\; \ldots \:$, 
$(\sigma_1^A \times B_M, \Sigma_1^A \times B_M)$,
in this order. It is easy to check that each of the steps above is a legitimate collapse: When we remove $\sigma_1^A \times B_i$  all the faces of the type $\sigma_1^A \times \beta$ containing $\sigma_1^A \times B_i$ have already been removed, because in the list $B_1, \ldots, B_M$ the face $\beta$ appears before $B_i$. On the other hand, $\sigma_1^A$ is a \emph{free} face of $\Sigma_1^A$, thus no face of the type $\alpha \times B_i$ may contain $\sigma_1^A \times B_i$ other than $\Sigma_1^A \times B_i$.

Next, we consider the \emph{second} pair of faces  $(\sigma_2^A, \Sigma_2^A)$ that appears in the collapse of $A$ onto $C_A$ and we repeat the procedure above, and so on: In the end, the only faces left are those of $C_A \times B$.
\end{proof}

\begin{cor} \label{cor:cohen}
If $A$ is collapsible, then $A \times B$ collapses onto a copy of $B$. 
\end{cor}

Since the product of the dunce hat with a segment $\mathbb{I}$ is collapsible \cite{ZeemanDunceHat}, the collapsibility of both $A$ and $B$ \emph{strictly} implies the collapsibility of $A \times B$. 

Now, consider a $1$-sphere $S$ consisting of four edges. The $2$-complex $S \times S$ is a cubical torus; after the removal of a facet, it collapses onto the union of a meridian and a longitude of the torus. (Topologically, a punctured torus retracts to a bouquet of two circles.) This can be generalized as follows:

\begin{proposition} \label{prop:uffa}
Let $A$ and $B$ be two polytopal complexes. Let $\Delta_A$ (resp. $\Delta_B$) be a facet of $A$ (resp. $B$). 
If $A - \Delta_A$ collapses onto some complex $C_A$ and if $B - \Delta_B$ collapses onto some complex $C_B$ then
$(A \times B) - (\Delta_A \times \Delta_B)$ collapses onto 
$(A \times C_B)  \cup (C_A \times B)$.
\end{proposition}

\begin{proof} We start by forming three ordered lists of pairs of faces.
Let 
$\left( \sigma_1, \Sigma_1\right)$,  
$\; \ldots \:$, 
$\left( \sigma_U, \Sigma_U\right)$
be the list of the removed pairs of faces in the collapse of $A - \Delta_A$ onto $C_A$. (We assume that higher dimensional faces are collapsed first.) 
Analogously, let $\left( \gamma_1, \Gamma_1\right)$,  
$\; \ldots \:$, 
$\left( \gamma_V, \Gamma_V\right)$
be the list of all the removed pairs in the collapse of $B - \Delta_B$ onto $C_B$.  Let then $B_1, \ldots, B_W$ be the list of all the faces of $B$ that are not in $C_B$, ordered by weakly decreasing dimension.

The desired collapsing sequence for $(A \times B) - (\Delta_A \times \Delta_B)$ consists of $U+1$  distinct phases:

\smallskip

\begin{compactdesc}
\item{$\!\!\!$\textsc{Phase} $0$:} We remove from $(A \times B) - (\Delta_A \times \Delta_B)$ the $V$ pairs of faces
$\left( \Delta_A \times\gamma_1, \Delta_A \times\Gamma_1\right)$,  
$\left( \Delta_A \times \gamma_2, \Delta_A \times \Gamma_2\right)$,
$\; \ldots, \:$,
$\left( \Delta_A \times\gamma_V, \Delta_A \times\Gamma_V\right)$,
  in this order. Analogously to the proof of Proposition \ref{prop:collapseproduct0}, one sees that all these removals are elementary collapses. They wipe away the ``$\Delta_A$-layer'' of $A \times B$, but not entirely: The faces $\alpha \times \beta$ with $\beta$ in $C_B$ are still present. What we have written is in fact a collapse of $(A \times B) - (\Delta_A \times \Delta_B)$ onto the complex 
$\left(\,(A - \Delta_A) \times B \right) \:\: \cup \:\: ( \Delta_A \times C_B )$.
\item{$\!\!\!$\textsc{Phase} $1$:} We take the first pair $\left( \sigma_1, \Sigma_1 \right)$ in the first list and we perform the $W$ elementary collapses  
$\left( \sigma_1 \times B_1, \Sigma_1 \times B_1\right)$,
$\; \ldots \:$, 
$\left( \sigma_1 \times B_W, \Sigma_1 \times B_W\right)$.
This way we remove (with the exception of $\Sigma_1 \times C_B$) the $\Sigma_1$-layer of $A \times B$, where $\Sigma_1$ is the first facet of $A$ to be collapsed away in $A - \Delta_A \searrow C_A$. 
\item{$\quad  \, \vdots$}
\item{$\!\!\!$\textsc{Phase} $j$:} We consider $\left( \sigma_j, \Sigma_j\right)$ and proceed as in Phase $1$, performing $W$ collapses to remove (with the exception of $\Sigma_j \times C_B$) the $\Sigma_j$-layer of $A \times B$.
\item{$\quad  \, \vdots$}
\item{$\!\!\!$\textsc{Phase} U:} We consider $\left( \sigma_U, \Sigma_U\right)$ and proceed as in Phase $1$, performing $W$ collapses to remove (with the exception of $\Sigma_U \times C_B$) the $\Sigma_U$-layer of $A \times B$.
\end{compactdesc}
Eventually, the only faces of $A \times B$ left are the polytopes of $A \times C_B \; \cup \; C_A \times B$.
\end{proof}

\begin{cor} \label{cor:collapseproduct} 
Given $s$ polytopal complexes $A_1, \ldots, A_s \;$, suppose that each $A_i$ after the removal of a facet collapses onto some lower-dimensional complex $C_i\;$. Then the complex $A_1 \times \ldots \times A_s$ after the removal of a facet collapses onto
\[ (C_1 \times A_2 \times \ldots \times A_s) 
\; \: \: \cup \: \: \;
(A_1 \times C_2 \times A_3 \times \ldots \times A_s)
\; \: \: \cup \; \: \: \ldots \; \: \: \cup \; \: \:
(A_1 \times \ldots \times A_{s-1} \times C_s) \;.
\]
In particular, if $\dim C_i = \dim A_i - 2$ for each $i \,$, then $A_1 \times \ldots \times A_s$ minus a facet collapses onto a complex of dimension 
$\dim A_1 + \ldots + \dim A_s - 2$.
\end{cor}

\begin{proof}
It follows from Proposition \ref{prop:uffa}, by induction on $s$.
\end{proof}

\begin{remark}
Proposition \ref{prop:collapseproduct0}, Proposition \ref{prop:uffa} and Corollary \ref{cor:collapseproduct} can be easily extended to the generality of finite regular CW complexes (see e.g. Bj\"{o}rner \cite[p.~1860]{BJOE} for the definition). 
\end{remark}

\begin{example} \label{ex:perelman}
  Let $C$ be the boundary of the three-dimensional cube $\mathbb{I}^3$; removing a square from $C$ one obtains a collapsible $2$-complex. The product $C \times C$ is a cubical $4$-manifold homeomorphic to $S^2 \times S^2$ (and not homeomorphic to $S^4$). 
The $4$-complex obtained by removing a facet from $C \times C$ collapses onto a $2$-complex, by Proposition \ref{prop:uffa}. Therefore, $C \times C$ is LC.  Note that the second homotopy group of $C \times C$ is nonzero. However, as observed by Durhuus and Jonsson \cite{DJ} \cite[Lemma~1.6.3]{Benedetti-diss}, every LC $d$-manifold is simply connected.
\end{example}

\begin{remark} \label{rem:perelman}
The previous example can be generalized by taking the product of the boundary of the $3$-cube $\mathbb{I}^3$ with the boundary of the $(d-1)$-cube $\mathbb{I}^{d-1}$ ($d \ge 4$). As a result, one obtains a cubical $d$-manifold that is homeomorphic to $S^2 \times S^{d-2}$ (and not homeomorphic to $S^d$). This $d$-ma\-nifold is LC, because the boundary of a $(d-1)$-cube is shellable and LC. In contrast, no manifold homeomorphic to  $S^1 \times S^{d-1}$ is LC, because LC manifolds are simply connected.
\end{remark}

\subsubsection*{Acknowledgements} The author wishes to thank G\"{u}nter Ziegler for helpful discussions.

\enlargethispage{2mm}

\end{document}